\documentclass[11pt]{amsart}
\usepackage{amssymb,amsmath,amsthm,amscd}
\usepackage[all]{xy}

\numberwithin{equation}{section}

\newtheoremstyle{fancy1}{10pt}{10pt}{\itshape}{12pt}{\textsc\bgroup}{.\egroup}{8pt}{}
\newtheoremstyle{fancy2}{10pt}{10pt}{}{12pt}{\itshape}{.}{8pt}{ }

\theoremstyle{fancy1}

\newtheorem{lem}[equation]{Lemma}
\newtheorem{prop}[equation]{Proposition}
\newtheorem{thm}[equation]{Theorem}

\newtheorem{main}{Theorem}
\newtheorem*{main*}{Theorem}

\newtheorem*{cor*}{Corollary}

\theoremstyle{fancy2}

\newtheorem*{rem*}{Remark}

\newcommand{\cref}[1]{Corollary~\ref{#1}}

\newcommand{\RP}{\mathbb{R\mkern1mu P}}
\newcommand{\CP}{\mathbb{C\mkern1mu P}}

\newcommand{\Sph}{\mathbb{S}}

\newcommand{\C}{{\mathbb{C}}}
\newcommand{\R}{{\mathbb{R}}}
\newcommand{\Z}{{\mathbb{Z}}}

\renewcommand{\H}{\ensuremath{\operatorname{H}}}

\newcommand{\E}{\ensuremath{\operatorname{E}}}
\newcommand{\F}{\ensuremath{\operatorname{F}}}
\newcommand{\G}{\ensuremath{\operatorname{G}}}

\newcommand{\SO}{\ensuremath{\operatorname{SO}}}

\renewcommand{\O}{\ensuremath{\operatorname{O}}}
\newcommand{\Sp}{\ensuremath{\operatorname{Sp}}}
\newcommand{\U}{\ensuremath{\operatorname{U}}}
\newcommand{\SU}{\ensuremath{\operatorname{SU}}}

\newcommand{\T}{\ensuremath{\operatorname{T}}}
\renewcommand{\S}{\ensuremath{\operatorname{S}}}

\newcommand{\N}{\ensuremath{\operatorname{N}}}
\newcommand{\K}{\ensuremath{\operatorname{K}}}
\renewcommand{\L}{\ensuremath{\operatorname{L}}}

\def\con#1=#2(#3){#1 \equiv #2 \bmod{#3}}

\newcommand{\diag}{\ensuremath{\operatorname{diag}}}

\newcommand{\codim}{\ensuremath{\operatorname{codim}}}

\newcommand{\Ad}{\ensuremath{\operatorname{Ad}}}

\DeclareMathOperator{\Iso}{Iso}

\DeclareMathOperator{\id}{id}

\newcommand{\no}{\noindent}
\newcommand{\co}{{cohomogeneity }}
\newcommand{\coo}{{cohomogeneity one }}
\newcommand{\com}{{cohomogeneity one manifold }}

\newcommand{\Kmo}{\K_{\scriptscriptstyle{0}}^{\scriptscriptstyle{-}}}
\newcommand{\Kpm}{\K^{\scriptscriptstyle{\pm}}}
\newcommand{\Kp}{\K^{\scriptscriptstyle{+}}}
\newcommand{\Km}{\K^{\scriptscriptstyle{-}}}

\newcommand{\subo}{_{\scriptscriptstyle{0}}}

\begin{document}

\newcommand{\spacing}[1]{\renewcommand{\baselinestretch}{#1}\large\normalsize}
\spacing{1.15}

\title{Symmetries of Eschenburg  spaces and the Chern problem}
\author[K.\ Grove, K.\ Shankar, W.\ Ziller]{Karsten Grove, Krishnan
Shankar, Wolfgang Ziller}

\dedicatory{Dedicated to the memory of S.\ S.\ Chern.}

\address{University of Maryland\\
        College Park , MD 20742}
\email{kng@math.umd.edu}

\address{University of Oklahoma\\
        Norman , OK 73019}
\email{shankar@math.ou.edu}

\address{University of Pennsylvania\\
        Philadelphia, PA 19104}
\email{wziller@math.upenn.edu}

\thanks{All three authors were supported by  grants from the National Science
Foundation and the third author  by the Francis J. Carey Term Chair
and the Clay Institute}

\maketitle



To advance our basic knowledge of manifolds with positive (sectional)
curvature it is essential to search for new examples, and to get a
deeper understanding of the known ones. Although any positively curved
manifold can be perturbed so as to have trivial isometry group, it is
natural to look for, and understand the most symmetric ones, as in the
case of homogeneous spaces. In addition to the \emph{compact rank one
symmetric spaces}, the complete list (see \cite{bergery}) of simply
connected homogeneous manifolds of positive curvature consists of the
\emph{Berger} spaces $B^7$ and $B^{13}$ \cite{berger}, the
\emph{Wallach} spaces $W^6, W^{12}$ and $W^{24}$ \cite{wallach}, and
the infinite class of so-called \emph{Aloff--Wallach} spaces,
$\mathcal{A}^7$ \cite{aw}. Their full isometry groups were determined
in \cite{shankar3}, and this knowledge provided new basic information
about possible fundamental groups of positively curved manifolds, and
in particular to counter-examples of the so-called \emph{Chern
conjecture} (see \cite{shankar} and \cite{grovshank, bazaikin2}).

Our purpose here is to begin a systematic analysis of the isometry
groups of the remaining known manifolds of positive curvature, i.e.,
of the so-called \emph{Eschenburg spaces}, $\mathcal{E}^7$
\cite{esch1,esch2} (plus one in dimension 6) and the \emph{Bazaikin
spaces}, $\mathcal{B}^{13}$ \cite{bazaikin}, with an emphasis on the
former.  Any member of $\mathcal{E}$ is a so-called \emph{bi-quotient}
of $\SU(3)$ by a circle:
$$
{\mathcal E}=\diag(z^{k_1}, z^{k_2},
z^{k_3})\backslash \SU(3)/ \diag(z^{l_1}, z^{l_2}, z^{l_3})^{-1}, \, |z|=1
$$
with $\sum k_i=\sum l_i$.  They contain the homogeneous Aloff--Wallach
spaces $\mathcal{A}$, corresponding to $l_i=0, i=1,2,3$, as a special
subfamily. Similarly, any member of $\mathcal{B}$ is a bi-quotient of
$\SU(5)$ by $\Sp(2)\S^1$ and the Berger space, $B^{13} \in
\mathcal{B}$. It was already noticed several years ago by the first
and last author, that both $\mathcal{E}$ and $\mathcal{B}$ contain an
infinite family $\mathcal{E}_1$ respectively $\mathcal{B}_1$ of
cohomogeneity one, i.e., their isometry groups have $1$-dimensional
orbit spaces (see section 1 and \cite{ziller}).  The work in
\cite{gwz} shows that up to diffeomorphism, these are the only
manifolds from $\mathcal{E}$ and $\mathcal{B}$ of cohomogeneity
one. There is a larger interesting subclass $\mathcal{E}_2
\subset\mathcal{E}$, corresponding to $l_1=l_2=0$, which contains
$\mathcal{E}_1$ as well as $\mathcal{A}$, and whose members have
cohomogeneity two.  We point out that $\mathcal{E}_1 \cap \mathcal{A}$
has only one member $A_{1,1}$, the unique Aloff--Wallach space that is
also a normal homogeneous space (see \cite{wilking:normal}).

From our results about isometry groups we get in particular

\begin{main}
The isometry group of any $E \in \mathcal{E}_2$ has dimension 11, 9, 7 or 5,
corresponding to the cases $E = A_{1,1}$, $E \in \mathcal{A} - \{A_{1,1}\} $,
$E \in  \mathcal{E}_1 - \{A_{1,1}\}$, or $E \in \mathcal{E}_2 - (\mathcal{E}_1
\cup \mathcal{A})$ respectively.
\end{main}

\no There are at most two possible groups for each class, and in all
cases we know explicitly what they are (at least up to components, see
Theorems \ref{isoE1}, \ref{isoB1} and \ref{isoE2}). We do not know if
there are any $E \in \mathcal{E} - \mathcal{E}_2$ that have
cohomogeneity two, but they all have a cohomogeneity four action. For
more information about equivalences up to diffeomorphism,
homeomorphism and homotopy within the classes $\mathcal{E}$ and
$\mathcal{B}$, we refer to \cite{shankar4}, \cite{escheretal} and
\cite{zilleretal}.  We note in particular that there are spaces from
the disjoint classes $\mathcal{E} - \mathcal{E}_2$ and $\mathcal{E}_2$
that are diffeomorphic, but our results are about isometry types.

The explicit description of each $E \in \mathcal{E}$ easily gives
rise to a ``natural group" of isometries $\N(E) \subset \Iso(E)$
(see section 1) of cohomogeneity $0, 1, 2$ or $4$. In the first case
we are done by \cite{shankar3}, and in the second by \cite{gwz},
after we have established what the so-called cohomogeneity one
diagrams are in this case. For the third case we show that if the
given action extends to a cohomogeneity 1 or 0 action, then $E$ is
already one of the more special spaces listed in $\mathcal{E}_1$ or
in $\mathcal{A}$. The general case of cohomogeneity four is not
treated in this paper.

\smallskip

Our concrete knowledge of the isometry groups of spaces $E \in
\mathcal{E}$ allows us to expand the list of positively curved
manifolds with interesting fundamental groups. Here subgroups of
$\SO(3)$ are particularly interesting, since many of them do not
occur as \emph{space form groups}, i.e., as fundamental groups of
spaces of constant curvature. The non-abelian simple group $A_5$ and
the abelian non-cyclic group $\Z_2 \oplus \Z_2$ are examples of such
groups. We will see that $\SO(3)$ itself acts freely and
isometrically on only one Aloff--Wallach space, and one Eschenburg
space (already found in \cite{shankar}). Nevertheless we will show
the following, which adds infinitely many spaces with distinct
homotopy types that violate Chern's conjecture for fundamental
groups of positively curved manifolds.

\begin{main}
For any finite subgroup $\Gamma \subseteq \SO(3)$, there exist
infinitely many spaces in $\mathcal{E}_1$ as well as in $\mathcal{E}_2
- \mathcal{E}_1$ on which $\Gamma$ acts freely and isometrically.

Moreover, for any odd positive integers $p$ and $q$ with gcd$(p+1,q) =
1$ the group $\Z_2 \times \Z_{2q}$ acts freely and isometrically on
$E_p \in \mathcal{E}_1$.
\end{main}

\smallskip

We have divided the paper into four sections. In the first section we
set up notation, including the precise definitions of the objects we
are interested in, and present the tools needed for our
proofs. Section 2 deals with the cohomogeneity one spaces
$\mathcal{E}_1$ and $ \mathcal{B}_1$. The bulk of our work is in
section 3 which provides a detailed analysis of the class
$\mathcal{E}_2$ (we also include a brief discussion for the
6-dimensional ``Eschenburg flag"). In section 4 we use our
knowledge of isometry groups developed in sections 2 and 3 to find
free isometric actions on manifolds from $ \mathcal{E}_2$.

This work was completed while the third author was visiting IMPA and
he would like to thank the Institute for its hospitality.


\section{Preliminaries and tools}

The general strategy for determining the isometry groups of the
Eschenburg and Bazaikin spaces has two steps. The first and fairly
simple step is to exhibit a (connected) group of isometries which
arises naturally from the description of the space. In the second and
more difficult step we then show that it cannot be enlarged to a group
of larger dimension. This uses an analysis of orbit spaces and
isotropy groups. The possible choices of enlargements are severely
restricted by a number of classification theorems about positively
curved manifolds with large isometry groups.  In this section we will
describe the first step, and provide the general tools needed for the
second step.  \bigskip

\begin{center}
\emph{Biquotient metrics and Natural Isometries.}
\end{center}

Throughout the paper, we let $\Iso(M)$ denote the full group of
isometries of a Riemannian manifold $M$. As usual, the identity
component of a Lie group $\G$ will be denoted by $\G_{\subo}$, and if
$\H \subset \G$ is a closed subgroup, then $\N_{\G}(\H)$ is the
normalizer of $\H$ in $\G$, or just $\N(\H)$ if it is clear from the
context.

By definition, a \textit{biquotient} manifold $M$ is the orbit space
$\G/\!/\U$ of a compact Lie group $\G$, by a subgroup $\U \subset
\G\times \G$ acting freely as
\begin{align*}
   \U \times \G &\rightarrow \G,   \quad   \quad     (u_1, u_2)\cdot g
   \rightarrow u_1
   \cdot g \cdot u_2^{-1},
\end{align*}
When $\U$ lies strictly in one factor of $\G\times \G$, then the
quotient is a homogeneous space.

The Riemannian metrics we consider on a biquotient $M = \G/\!/\U$
are always induced from a left invariant, $\Ad(\K)$-invariant metric on
$\G$ where  $\U \subset \G\times \K$, and $\K\subset \G$ is a closed
subgroup. We then have the inclusions $\Iso(M) \supset
\N_{\Iso(\G)}(\U)/\U \supset \N_{\G\times \K}(\U)/\U$ since $
\G\times \K \subset \Iso(\G)$. We will refer to $\N(M) :=
\N_{\G\times \K}(\U)/\U$ as the \emph{natural group of isometries}
of the biquotient $M = \G/\!/\U$.

\bigskip

\begin{center}
\emph{Eschenburg and Bazaikin spaces.}
\end{center}

We will now describe the special biquotients we are dealing with in
this paper, namely the Eschenburg and the Bazaikin spaces
$\mathcal{E}$ and $\mathcal{B}$. For both classes $\G = \SU(n)$ and
$\K = \U(n-1) = \S(\U(n-1)\U(1))$, where $n = 3$ and $5$ respectively.
From the above discussion this will already determine the metrics we
consider on the orbit spaces $\G/\!/\U$, for $\U \subset \G \times
\K$.

To describe the spaces in $\mathcal{E}^7$, we proceed as follows:

\no Let $\bar{a}:=(a_1,a_2,a_3)$, $\bar{b}:=(b_1,b_2,b_3)$ be triples of
integers such that $\sum a_i = \sum b_i : = c$. Let
$$
   \S^1_{\bar{a}, \bar{b}} =
\{(\diag(z^{a_1},z^{a_2},z^{a_3}), \diag(z^{b_1},z^{b_2},z^{b_3}))
\mid z \in \U(1) \}
$$

\noindent The $\S^1_{\bar{a}, \bar{b}}$ action on $\SU(3)$ is free if
and only if for every permutation $\sigma \in \text{S}_3$, $\gcd(a_1 -
b_{\sigma(1)}, a_2 - b_{\sigma(2)}) = 1$. In this case, we will call
the resulting 7-manifold, $E_{\bar{a},\bar{b}} := \SU(3) /\!/
\S^1_{\bar{a},\bar{b}}$, an \textit{Eschenburg space}. Note that
$\S^1_{\bar{a}, \bar{b}} \subsetneq \SU(3) \times \U(2)$, but its
action is the same, up to an ineffective kernel, as that by
$\S^1_{3\bar{a} - \bar{c}, 3\bar{b} - \bar{c}} \subset \SU(3) \times
\U(2)$, where $\bar{c} = (c,c,c)$. In \cite{esch1} it was shown that
the Eschenburg metric on $E_{\bar{a},\bar{b}}$ has positive sectional
curvature if and only if one of the following holds:
\begin{equation}\label{pos}
b_i \not\in [a_{\min}, a_{\text{max}}], \text { or } a_i
\not\in [b_{\text{min}}, b_{\text{max}}] \text{ for all } i.
\end{equation}

\no Strictly speaking, we need to allow the invariance of the metric
to be switched, and to choose any of the 3 different block embeddings
of $\U(2)\subset\SU(3)$ in order to obtain this necessary and
sufficient condition. But for convenience, we will fix the embedding
and assume the metric is left invariant. We reserve the notation
$\mathcal{E}$ for those Eschenburg spaces $E_{\bar{a},\bar{b}}$ that
have positive curvature.  If the action by $\S^1_{\bar{a}, \bar{b}}$
is only one sided, we obtain the subfamily of homogeneous
Aloff--Wallach spaces $A_{k,l}=\SU(3)/\diag(k,l,-(k+l))$ with
$\gcd(k,l)=1$. Here we can assume, up to conjugacy and change of
orientation, that $k\ge l\ge 0$. $A_{k,l}$ admits a homogeneous metric
with positive curvature if and only if $l>0$.

\bigskip

To describe the spaces in $\mathcal{B}$, consider a five tuple of
integers $\bar{p} = (p_1, p_2, p_3, p_4, p_5)$ with $q: = \sum
p_i$. Let
$$\Sp(2)\S^1_{\bar{p}} =
\{(\diag(z^{p_1},z^{p_2},z^{p_3},z^{p_4},z^{p_5}), \diag(\Sp(2), z^{q})\},$$

\no where $\Sp(2) \subset \SU(4)$ is embedded in the upper block of
$\SU(5)$. The action of $\Sp(2)\S^1_{\bar{p}}$ on $\SU(5)$ is free if
and only if all $p_i$ are odd and for all permutations $\sigma \in
{\rm S}_5$, $\gcd(p_{\sigma(1)} + p_{\sigma(2)}, p_{\sigma(3)} +
p_{\sigma(4)}) = 2$. In this case, we say that $B_{\bar{p}}:=
\SU(5)/\!/\Sp(2)\S^1_{\bar{p}}$ is a \textit{Bazaikin space}. As for
the Eschenburg spaces above we note that $\Sp(2)\S^1_{\bar{p}}
\subsetneq \SU(5)\times \U(4)$, but its action is the same as that of
$\Sp(2)\S^1_{5\bar{p} - \bar{q}} \subset \SU(5)\times \U(4)$, where
$\bar{q} = (q,q,q,q,q)$. From the treatment in \cite{ziller} of
Bazaikin's work \cite{bazaikin}, we know that the Eschenburg metric on
$B_{\bar{p}}$ has positive curvature if and only if
\begin{equation}\label{Bapos}
 p_{\sigma(1)} + p_{\sigma(2)} >0  \text{ for all permutations }
 \sigma \in {\rm S}_5.
\end{equation}

\no We reserve the notation $\mathcal{B}$ for those Bazaikin spaces
$B_{\bar{p}}$ that have positive curvature. In the case of
$\bar{p}=(1,\dots ,1)$ we obtain the unique Bazaikin space which is
homogeneous, the Berger space $B^{13}=\SU(5)/\Sp(2)\S^1$.

\bigskip

\begin{center}
\emph{Group Enlargements.}
\end{center}

In this subsection we consider the situation where an isometric $\G$
action on $M$ is a sub-action of an isometric $\G^*$ action, and $\G
\subseteq \G^*$ and $M$ are all compact and connected. Clearly then,
one has an induced submetry $\pi: M/\G \longrightarrow M/\G^*$ and
$\dim(M/\G^*) \le \dim(M/\G)$.  Moreover, if we let $(M/\G)_0$ denote
the \emph{regular} part of $M/\G$, corresponding to the \emph{principal}
$\G$ orbits $M_0$ in $M$, and similarly for the $\G^*$ action, we
have:

\begin{lem}[Submetry]\label{submetry}
All principal $\G^*$ orbits in $M$ are equivalent as $\G$ manifolds as
well. Moreover, the subset
$(M/\G)_0 \cap \pi^{-1}(M/\G^*)_0$ is open and dense in $M/\G$, and the image
$\pi((M/\G)_0 \cap \pi^{-1}(M/\G^*)_0) =
(M/\G^*)_0$. In particular, $M/\G$ and $M/\G^*$ are isometric if $\dim(M/\G) =
\dim(M/\G^*)$.
\end{lem}

\begin{proof}
Let $P$ and $P^*$ denote the projections from $M$ to $M/\G$ and $M/\G^*$
respectively. We have the following commutative diagram.
$$
  \xymatrix{
    M \ar@{->}[r]^{P} \ar@{->}[dr]_{P^*} & M/\G \ar@{->}[d]^{\pi} \\
     & M/\G^*
     }
$$

The collection $M^*_0$ of all principal $\G^*$ orbits in $M$ is an
open and dense $\G$ invariant subset of $M$. Moreover, all principal
$\G^*$ orbits are equivalent as $\G$ manifolds as well. Clearly then,
the set of principal $\G$ orbits in $ M^*_0$ is open and dense in $M$,
and in fact $P^*(M_0 \cap M^*_0) = (M/\G^*)_0$ and $P(M_0 \cap M^*_0)
= (M/\G)_0 \cap \pi^{- 1}(M/\G^*)_0$.

Now suppose $\dim(M/\G) = \dim(M/\G^*)$. Then the map $\pi$ is a local
isometry from $ (M/\G)_0 \cap \pi^{-1}(M/\G^*)_0$ onto
$(M/\G^*)_0$. Since $\G$ and $\G^*$ are connected, it is also clearly
1-1, and hence an isometry. It now follows from the first part that
$M/\G$ and $M/\G^*$ are isometric under $\pi$.
\end{proof}

\bigskip

\bigskip

\begin{center}
\emph{Size Restrictions.}
\end{center}

In the presence of positive curvature, the size of the isometry group
is restricted, which will be an important tool in our discussions.

Recall that the \emph{symmetry rank} of a manifold by definition is
the rank of its isometry group. In positive curvature this rank is
bounded above by \cite{grove-searle:rank}.

\begin{thm}[Rank Rigidity]\label{rank}
Assume that  a $k$-dimensional torus acts effectively and
isometrically on a positively curved simply connected $n$-manifold
$M$. Then $k \le [\frac{n+1}{2}]$, and equality holds only when $M$ is
diffeomorphic to $\Sph^n$ or $\CP^{n/2}$.
\end{thm}

Another measurement for the size of a group is its dimension. The
\emph{degree of symmetry} of a Riemannian manifold $M$ is by
definition the dimension of its isometry group. This dimension is
severely restricted in positive curvature by the following result of
Wilking \cite{wilking3}

\begin{thm}[Symmetry Degree]\label{degree}
Let $(M^n, g)$ be a simply connected, Riemannian
manifold of positive curvature. If the symmetry degree of $M^n$ is at
least $2n-5$, then $M^n$ is homotopy equivalent to a compact, rank one
symmetric space, or $M^n$ is isometric to a homogeneous space of
positive sectional curvature.
\end{thm}

The \emph{cohomogeneity}, i.e., the dimension of the orbit space gives
yet another measurement for the size of a transformation group. A
related invariant is the so-called \emph{fixed point cohomogeneity}
which is the dimension of the normal sphere to the fixed point set in
the orbit space. A manifold that supports an action of fixed point
cohomogeneity 0 is called \emph{fixed point homogeneous}. Although
this will not be used in the sense of size here, the following
classification result of \cite{grove-searle:cofix} is quite useful for
our investigations:

\begin{thm}[Fixedpoint Homogeneity]\label{cofix}
A fixed point homogeneous simply connected  manifold of positive curvature is
diffeomorphic to a rank one symmetric space.
\end{thm}


\section{Cohomogeneity One}

In this section we single out the subclasses $\mathcal{E}_1 \subset
\mathcal{E}$ and $\mathcal{B}_1 \subset \mathcal{B}$ of positively
curved cohomogeneity one Eschenburg and Bazaikin spaces, and determine
their full isometry groups.  Here
\begin{equation}
\mathcal{E}_1 = \{E_p = E_{\bar{a},\bar{b}} \in \mathcal{E} \mid
\bar{a} = (1,1,p), \  \bar{b} = (0,0,p+2), \   p> 0\}
\end{equation}

\no and
\begin{equation}
\mathcal{B}_1 = \{B_p = B_{\bar{p}} \in \mathcal{B}  \mid  \bar{p} =
(1,1,1,1,2p-1),  \   p> 0\}
\end{equation}

\no From \eqref{pos} and \eqref{Bapos}  we know that all these
manifolds have positive curvature when equipped with the Eschenburg
biquotient metric. We also point out that the $\S^1_p$, and
$\S^1_{-p-1}$ actions on  $\SU(3)$ are equivalent via the inverse map
of $\SU(3)$. Moreover, $E_0 \approx E_{-1}$ only has non-negative
curvature in the Eschenburg metric, and in fact does not support any
cohomogeneity one metric of positive curvature by \cite{gwz}. Note
also, that $E_1$ is the homogeneous Aloff--Wallach space $A_{1,1}$ and
$B_1$ is the homogeneous Berger space $B^{13}$.

To see that each $E_p$ has cohomogeneity one note that the natural
action by $\U(2)\times \SU(2)$ (as well as by $\SU(2)\times \U(2)$)
on $\SU(3)$ commutes with the $\S^1_p$ action, and that $\U(2)
\backslash \SU(3)/ \SU(2)$ = $\CP^2/\SU(2)$, which is an interval.
We also note that  $\SU(2)\times\SU(2)$ and $\S^1_p$ generates
$\U(2)\times \SU(2)$ and hence the sub-action by
$\SU(2)\times\SU(2)$ is  cohomogeneity one as well.

In the case of $B_p$ we see that the natural left action by $\U(4)$ on
$\SU(5)$ commutes with the action of $\Sp(2)\S^1_p$ which in turn has
cohomogeneity one. Indeed, the induced action by $\Sp(2)$ on
$\U(4)\backslash \SU(5) = \CP^4$ is the standard sub-action of
$\Sp(2)\subset \SU(4)\subset \SU(5)$ and hence, $\U(4)\backslash
\SU(5)/ \Sp(2)$ is an interval. As in the case of Eschenburg spaces,
we note that the group generated by $\SU(4)$ and $\Sp(2)\S^1_p$ is the
same as the one generated by $\U(4)$ and $\Sp(2)$, and therefore the
sub-action by $\SU(4)$ is also of cohomogeneity one.

\smallskip

Using the tools from section 1 one can prove that if any of these
actions extend to a transitive isometric action then $p = 1$ for both
classes. It is actually known that none of $E_p$ or $B_p$ for $p >1$
is even homeomorphic to a homogeneous space (see \cite{shankar4},
\cite{zilleretal}).  Since the full isometry groups of these
homogeneous spaces were determined in \cite{shankar} it remains to
consider $E_p$ and $B_p$ for $p >1$.

We first determine the identity component of $\Iso(E_p)$ by
analyzing the sub-action by $\G=\SU(2) \times \SU(2)\subset
\U(2)\times \SU(2)$ , which we noted above  is also  cohomogeneity
one. It is important for us, however, to determine the associated
\emph{group diagram}, $\H \subset \{\Km,\Kp\} \subset \G$, i.e., the
isotropy groups along a minimal geodesic between the two
non-principal orbits $B_{\pm} = \G/\Kpm$ corresponding to the end
points of the orbit space interval. This information is also used as
a \emph{recognition tool} in the classification work of \cite{gwz}.

\begin{prop}
The cohomogeneity one action of  $\G = \SU(2)\times \SU(2)$ on $E_p$
has principal isotropy group
$\H=\{(\pm{\id})^{p+1},(\pm{\id})^{p}\}\cong \Z_2$ and singular
isotropy groups $\Km = \Delta \SU(2)\cdot\H$ and $\Kp =
\S^1_{(p+1,p)}$ embedded with slope  $(p+1, p)$ in a maximal torus
of $\SU(2)\times \SU(2)$.
\end{prop}

\begin{proof}
Consider  the point $p_- = \S^1_p(e) = \{\diag(z, z, \bar{z}^2)\}$  in
$E_p$, and let $B_- = \G(p_-)$ be the orbit of this point under the
action of $\SU(2)\times \SU(2)$. The identity component of the
isotropy group at $p_-$ is clearly $\Kmo = \Delta \SU(2)$ and the
second component in $\K$ is generated by $(\id,-\id)$. Hence $B_-
\cong \SO(3) \cong \RP^3$.

Since we already saw that $E_p$ is cohomogeneity one, the action of
$\Km$ (effectively by $\Kmo=\SU(2)$) on the normal space
$T_{p_-}^{\perp}$ to $B_-$ at $p_-$ is the standard action of $\SU(2)$
on $\C^2 \cong \R^4$.  Because this action of $\SU(2)$ restricted to
the normal 3-sphere is both transitive and free, the effective version
of the action by $\G$ has trivial principal isotropy group.

To find the other singular orbit $B_+$, note that $v_- = E_{13} \in
\mathfrak{su(3)}$ (standard basis element for the skew symmetric
matrices) represents a normal vector to $B_-$ at $p_-$. One easily
checks that the one parameter group $\exp(tv_-)$ is still a geodesic
in the left invariant Eschenburg metric on $\SU(3)$ and hence on $E_p$
as well. It intersects $B_-$ again at time $\pi$, and not
earlier. This implies that $p_+ = \exp(\frac{\pi}{2}v_-) \in B_+$, and
$p_+$ is represented by $E_{13} + \diag(0,1,0)$.  To determine $\Kp$,
let $(g_1,g_2) \in \SU(2)\times \SU(2) \subset \SU(3)\times\SU(3)$. We
identify $\SU(2)$ with the unit sphere in $\C^3$ as usual and let
$(a,b)$ correspond to $g_1$ and $(\alpha, \beta)$ to $g_2$.  Then
$(g_1, g_2) \in \Kp$ if and only if $(g_1, g_2)\cdot p_+ \in
\S^1_p(p_+)$.  This implies that $b = \beta =0$, $a = \bar{z}^{p+1}$
and $\alpha = \bar{z}^p$. Or equivalently, $a = z^{p+1}$ and $\alpha =
z^p$.  Notice that if $p$ is even, then $\S^1_{(p+1, p)}$ goes through
$(-\id, \id)$ while it goes through $(\id, -\id)$ if $p$ is odd and
hence $\H=\{(\pm{\id})^{p+1},(\pm{\id})^{p}\}$.
\end{proof}

The action of $\G$ is ineffective with kernel $\H$, and hence the
natural group of isometries is $\U(2)\times \SO(3)$ when $p$ is odd,
and $\SO(3)\times \U(2)$ when $p$ is even. Furthermore, complex
conjugation on $\SU(3)$ normalizes the circle action and hence
induces an isometry as well.

\begin{thm}\label{isoE1}
For any integer $p > 1$, the full isometry group of $E_p$ is given by
$$
   \Iso(E_p) = (\U(2) \rtimes \Z_2)\times \SO(3),
$$
where  the second component is induced by complex conjugation in $\SU(3)$.
 \end{thm}

\begin{proof}
We already saw that $\U(2)\times \SO(3) \subset \Iso_0(E_p)$. Since
$E_p,\, p>1$ is not diffeomorphic to a homogeneous space
(\cite{shankar4}), any extension of the group action must again be
of \co one. Moreover, since the two singular orbits $B_-$ and $B_+$ are
non-isometric (not even of the same dimension), it must act trivially
on the orbit space interval. In other words any extension of the action
will have the same orbits. We now consider the codimension two orbit
$B_+=\S^3\times \S^3/\S^1_{(p+1,p)}$ and claim that the action of $
\Iso(E_p)$ on $B_+$ must be effective. Indeed, assume that an element
$\gamma\in \Iso(E_p)$ acts trivially on $B_+$. Then  $B_+$ is either
totally geodesic, or $\gamma$ acts by reflection on the two
dimensional normal space to $B_+$. The former is impossible since
$B_+$ does not support a homogeneous metric of positive curvature.
The latter is impossible as well, since $E_p$ would then have a
totally geodesic hypersurface, but in positive curvature this is only
possible when the manifold is homeomorphic to the sphere or to real
projective space.

We will now examine the size of the isometry group for a homogeneous
metric on $B_+$. We already know that $\G=\S^3\times\S^3\times\S^1$
(effectively  $\U(2)\times \SO(3)$) acts by isometries on $B_+$.
From the classification of 5-dimensional homogeneous manifolds it
follows that  it is also the identity component of the isometry
group, i.e., no larger connected group can act transitively. One
easily checks that  $\L=\{\diag(z^{p+1},z^p,1), \diag(1,w,w)\}$ is
the full isotropy group of $\S^3\times\S^3\times\S^1$ and that
$\N(\L)/(\L\cdot Z(\G))$ is trivial. Furthermore,   there exists an
outer automorphism of $\S^3\times\S^3\times\S^1$, unique up to inner
automorphisms, which preserves $\L$. From \cite[Theorem 3.1]{WZ},
(cf. also \cite{shankar3}), it then follows that the isometry group
of any homogeneous metric on $B_+$ can have at most two components.
Altogether, this completes the proof.
\end{proof}

\smallskip

We now turn to the isometry groups of the Bazaikin spaces $B_p \in
\mathcal{B}_1 , \, p>1$.  We already saw that the natural left $\U(4)$
action on $\SU(5)$ induces a cohomogeneity one action on
$B_p$. Furthermore, complex conjugation in $\SU(5)$ generates a second
component.

\begin{thm}\label{isoB1}
 The full isometry group of the cohomogeneity one manifold $B_p, \,
 p>1$ may be written as:
$$
   \Iso(B_p) = \U(4) \rtimes \Z_2,
$$
\no where the second component is induced by complex conjugation in
$\SU(5)$.
\end{thm}

\begin{proof}
We proceed as in the case of Eschenburg spaces above. Since $B_p,\,
p>1$ is not diffeomorphic to a homogeneous space (\cite{zilleretal}),
any extension of the group action will have cohomogeneity one and in
fact the same orbits when the singular orbits are different.

To determine the orbit structure of the action (cf. also \cite{gwz})
we consider the orbit equivalent sub-action by
$\G=\SU(4)\subset\U(4)$. The orbit through the identity is $
\SU(4)/(\Sp(2)\cup i\Sp(2))=\RP^5$ and the action by $\Km = \Sp(2)\cup
i\Sp(2)$ on the slice is nontrivial and hence
$\H=\SU(2)\cdot\Z_2$. Since $B_p$ is simply connected, the \coo action
cannot have any exceptional orbits and since $\Kp/\H$ is a sphere and
$\H$ is not connected, it must be one dimensional, i.e., $\Kp =
\SU(2)\cdot\S^1$. Since the centralizer of $\H$ in $\G$ is two
dimensional, $\S^1$ is allowed to have slopes $(r,s)$ inside this
two-torus. These slopes are then determined from those for the
Eschenburg spaces,since $E_p$ is the fixed point set of   the
involution $\diag(1,1,-1,-1,1)\in\SU(5)$, assuming that
$\SU(2)\subset\H$ is the lower $2\times 2$ block of $\SU(5)$,
(see \cite{taimanov}). It thus follows that $(r,s)=(p+1,p)$. Note
that since $E_p$ is not homogeneous, this observation also provides
a simple geometric proof that $B_p$ cannot be homogeneous, since
a totally geodesic submanifold of a homogeneous space is itself homogeneous
(cf.\ \cite{kobnom2}, Chapter VII, Corollary 8.10).

We again have that $\Iso(B_p)$ acts effectively on $B_+=\G/\Kp$ since
$B_+$ does not support a homogeneous metric with positive curvature,
and next determine the isometry group of the metric on $B_+$. To see
that $\U(4)$ is the identity component, one uses \cite[Theorem
4.1]{On} to show that $\U(4)$ cannot be enlarged to a bigger
transitive action. By computing the isotropy representation of
$\G/\Kp$, it follows that $\N(\Kp)\subset \SU(2)\cdot\U(2)$ and hence
$\N(\Kp)/\Kp=\S^1$ is connected. Since furthermore, $\U(4)$ has, up to
inner automorphisms, a unique outer automorphism, \cite[Theorem
3.1]{WZ} finishes the proof.
\end{proof}

\section{Cohomogeneity Two}

In this section we define and analyze a subclass $\mathcal{E}_2
\subset \mathcal{E}$ of positively curved, cohomogeneity two Eschenburg
spaces, and determine their isometry groups.  We will also briefly
treat the single 6-dimensional Eschenburg space, $E^6$.

Define
\begin{equation}
\mathcal{E}_2 = \{E_{\bar{p}} = E_{\bar{a},\bar{b}} \in \mathcal{E}
\mid \bar{a} = (p_1,p_2,p_3), \  \bar{b} = (0,0, p_1+p_2+p_3)\}
\end{equation}

We note that the $\S^1_{\bar{p}}$ action on $\SU(3)$ is free if and
only if $\gcd(p_i,p_j)=1$ for all $i\ne j$. In particular at most one
$p_i$ is even.  From \eqref{pos} it easily follows that the Eschenburg
metric has positive curvature, if and only if, up to reordering of the
$p_i$'s and changing the sign of all three, one of the following
holds:
\begin{equation}\label{coh2pos}
0 < p_1 \leq p_2 \leq p_3\, \text{ or }\,    0 < p_2 \leq p_3 \text{ and } p_1 <-p_3
\end{equation}

This class obviously contains the cohomogeneity one Eschenburg spaces,
$E_p \in \mathcal{B}_1$, where $\bar{p} = (1,1,p)$. Moreover, it
contains the (homogeneous) Aloff--Wallach spaces $\mathcal{A}$, since
$A_{p_1,p_2} = E_{\bar{p}}$ when $p_3 = -(p_1 +p_2)$.

\smallskip

Note that the natural action on $\SU(3)$ by $\T^2 \times \U(2)$
commutes with the $\S^1_{\bar{p}}$ action, and that $\T^2 \backslash
\SU(3)/\U(2) = \T^2 \backslash \CP^2$ is a right angled
triangle. Since $\S^1_{\bar{p}}$ and $\T^2 \times \SU(2)$ generate
$\T^2 \times \U(2)$, we see that the induced action by $\T^2 \times
\SU(2)$ on $E_{\bar{p}}$ has cohomogeneity two.  Clearly the effective
group is $\T^2 \times \SO(3)$ when all $p_i$ are odd, and $\T^2 \times
\SU(2)$ otherwise. For convenience we will work with $\T^2 \times
\SU(2)$ directly.  Although we do not need the full orbit structure of
this \co two action, the following information will be crucial.

\begin{lem}\label{lensspaces}
The orbits of the $\T^2 \times \SU(2)$ action on $E_{\bar{p}}$,
corresponding to the vertices in the quotient triangle, are lens
spaces with fundamental groups of order $|p_i+p_j|$.
\end{lem}

\begin{proof}
Let us first consider the orbit $B_1$  going through the image of
$\id \in \SU(3)$ in $E_{\bar{p}}$. The element
$((w_1,w_2),\diag(r,\bar{r}))\in \T^2 \times \SU(2)$ lies in the
isotropy of this point if and only if
$\diag(w_1r,w_2\bar{r},\bar{w}_1\bar{w}_2) =
\diag(z^{p_1},z^{p_2},\bar{z}^{p_1+p_2})$ for some $z\in\U(1)$. Since
$w_1$ and $w_2$ can be described arbitrarily, the subgroup
$1\times\SU(2)$ acts transitively on $B_1$ and has isotropy group
$\diag(z^{p_1},z^{p_2})$ with $z^{p_1+p_2}=1$, i.e., $B_1$ is a lens
space with fundamental group $\Z_{|p_1+p_2|}$.

The two other vertices can be represented by $E_{13} + \diag(0,1,0)$
and $E_{23} + \diag(1,0,0)$ since they are the endpoints of a geodesic
of length $\pi/2$ with initial vectors orthogonal to $B_1$ and
singular with respect to the isotropy action on the normal space of
$B_1$. The same computation as above shows that the orbits through
these two points are lens spaces with fundamental group
$\Z_{|p_1+p_3|}$ and  $\Z_{|p_2+p_3|}$ respectively.
\end{proof}

In order to determine the full isometry group, we need  to prove the
following claims:

\begin{itemize}
\item
If the $\T^2 \times \SU(2)$ action on $E_{\bar{p}}\in \mathcal{E}_2$
extends to an isometric  cohomogeneity one action, then $p_1 = p_2 =
0$, i.e., $E_{\bar{p}} \in \mathcal{E}_1$.
\item
If the $\T^2 \times \SU(2)$ action on $E_{\bar{p}}\in \mathcal{E}_2$
extends to a transitive isometric  action, then $p_3 = -p_1 - p_2 $,
i.e., $E_{\bar{p}} \in \mathcal{A}$.
\end{itemize}

\no We point out that classification results for positively curved
manifolds in cohomogeneity zero and one immediately yield
diffeomorphism conclusions in the above cases,  but our results are
about equality of the integer parameters.

\smallskip

Recall that by the \emph{rank rigidity theorem} \eqref{rank},
$\Iso(E_{\bar{p}})$ must have rank 3, and by the \emph{degree theorem}
\eqref{degree}, $\dim \Iso(E_{\bar{p}}) \le 8$, unless $E_{\bar{p}}$
is isometric to a homogeneous space with positive curvature. This
leaves only the following possible (almost effective) connected
extensions $\G^*$ of $\G = \T^2\times\SU(2)$:

\begin{itemize}
\item
$\G^* = \S^1\times \SU(2)\times\SU(2)$
\item
$\G^*= \S^1\times\SU(3)$, or $\G^* = \SU(2)\times\SU(3)$

\end{itemize}

We first deal with the latter \emph{extension case}:

\begin{prop}[Homogeneous Case]
Let $E_{\bar{p}} \in \mathcal{E}_2$, and assume the action of $\G =
\T^2 \times \SU(2)$ extends to a larger (connected) group $\G^*$ with
$\dim(\G^*) > 8$, then one of the following must occur:
\begin{itemize}
\item[(\textit{i})] $\G^* = \SU(2) \times \SU(3)$, $\bar{p} =
(1,1,-2)$ and $E_{\bar{p}}$ is the Aloff--Wallach space $A_{1,1}$.

\item[(\textit{ii})] $\G^* = \S^1 \times \SU(3)$, $\bar{p} = (p_1,
p_2, -p_1 -p_2)$ and $E_{\bar{p}}$ is the Aloff--Wallach space
$A_{p_1,p_2}$.
\end{itemize}
\end{prop}

\begin{proof}
The degree theorem implies that under the assumption $\dim(G^*) > 8$,
the Eschenburg space must be isometric to a homogeneous space. Since
among the positively curved homogeneous spaces in dimension $7$, only
the Aloff--Wallach spaces have possibly the same homotopy type as an
Eschenburg space (\cite{esch1}), it only remains to check the claims
about $\bar{p}$.

In the first case we note that only $A_{1,1}$ has an 11-dimensional
isometry group (cf.\ \cite{shankar3}). To see that indeed $\bar{p} =
(1,1,-2)$, we use the fact that the fourth cohomology group of
$E_{\bar{p}}$ is a finite cyclic group of order
$r = |p_1p_2+p_1p_3+p_2p_3|$ (cf.\ \cite{esch3}). Using the positive
curvature condition \eqref{coh2pos}, one easily sees that $r=3$ is
only assumed in the case of $\bar{p} = (1,1,-2)$.

In the second case, observe that there is, up to conjugacy, only one
immersed subgroup $\T^2\times \SU(2)$ in $\S^1\times\SU(3)$. Thus the
\co two action agrees with the one considered above on the
Aloff--Wallach spaces $A_{k,l}$ and hence the vertices in the orbit
space correspond to 3-dimensional lens spaces with fundamental groups
of orders $|k|, |l|, |k+l|$. On the other hand, for the action of
$\T^2\times \SU(2)$ on $E_{\bar{p}}$ the vertices correspond to lens
spaces whose fundamental groups have orders $|p_i+p_j|$. This imposes
severe restrictions on the $p_i$'s, and one easily shows that under
the positive curvature condition \eqref{coh2pos} this is only possible
when $\sum p_i = 0$.
\end{proof}

For the first \emph{extension case}, we will show:

\begin{prop}[Cohomogeneity One Case]
Let $E_{\bar{p}} \in \mathcal{E}_2$, and assume that the (almost
effective) action of $\G = \T^2 \times \SU(2)$ extends to a larger
(connected) group $\G^*$ with $ \dim(\G^*) \le 7$. Then $\G^* = \S^1
\times \SU(2)\times \SU(2)$ and $\bar{p} = (1,1,p) , p>0$.
\end{prop}

\begin{proof}
We have already seen that $\G^* = \S^1 \times \SU(2)\times \SU(2)$. It
remains to show that the action by $\G^*$ must be of cohomogeneity
one, and then to recognize $\bar{p}$.

If the $\G^*$ action also has cohomogeneity two, we see from the
\emph{submetry lemma} \eqref{submetry}, that it is orbit equivalent to
the $\G$ action. Since the latter must have finite principal isotropy
group, the fundamental group of the common principal orbits have
$\Z$-rank at least 2. On the other hand, any homogeneous quotient of
$\G^*$ has fundamental group with $\Z$-rank at most 1, and we conclude
that the $\G^*$ cation must have cohomogeneity one.

Now consider the submetry $\pi: M/\G \to M/\G^*$. From the submetry
lemma, we see that the inverse image of one of the end points of the
interval $M/\G^*$ is a vertex of the right angled triangle $M/\G$, and
that the opposite side is the inverse image of the other end
point. Let us denote the common singular $\G$ and $\G^*$ orbit by $B_1
= B_-^*$, and the other singular $\G^*$ orbit by $B_+^*$. Then $B_+^*$
is a cohomogeneity one $\G$ manifold at maximal distance to $B_1$.

As we saw, $B_1$ is a lens space with fundamental group of order
$|p_i+p_j|$ for some $i,j$ and now the bigger group $\G^*=\S^1 \times
\SU(2) \times \SU(2)$  also acts transitively on $B_1$. This implies
that $B_1$ is either $\Sph^3$ or $\RP^3$, i.e., that $|p_i+p_j| \le
2$. Indeed, if not, the $\G^*$ action on $B_1$ has a kernel $C$ of
dimension at least 3 since the isometry group of any lens space
$\Sph^3/\Z_{m}$, with $m \ge 3$ has dimension at most 4. In
particular, $C$ must contain an $\SU(2) \subset \G^*$. Since $B_1$ has
codimension 4, this $\SU(2)$ acts either trivially or transitively on
the normal spheres to $B_1$. In the first case $\SU(2)$ then acts
trivially on $M$, and in the second case $M$ is fixed point
homogeneous. Either one yields a contradiction (the second by
\ref{cofix}).

Our next claim is that the two other vertex $\G$ orbits, $B_2$ and
$B_3$, are lens spaces with isomorphic fundamental groups. To see
this, consider the cohomogeneity one $\G$ manifold $B_+^* =
\G^*/{\Kp}^*$. The singular orbits $B_- = B_2$ and $B_+ = B_3$ are
again 3-dimensional lens spaces with isotropy group $\T^2$, and the
principal orbits have isotropy $\S^1$. In particular $\dim B_+^* =
5$, and since $\codim B_-^* = 4 > 2$ it follows that $B_+^*$ is
simply connected by transversality . Thus ${\Kp}^* = \T^2$ and
$B_+^* = \S^1\times \SU(2)\times\SU(2)/\T^2 =
\SU(2)\times\SU(2)/\S^1_{s,t}$ for some $s,t$. Note also, that since
the $\SU(2) \subset \G$ commutes with $\T^2$ it must be one of the
$\SU(2)$ factors of $\G^*$, and we have a sub-action by $\G' =
\S^1\times\SU(2) \subset \SU(2)\times\SU(2)$ which is orbit
equivalent to the $\G$ action on $B_+^*$. The isotropy groups of the
$\G'$ action are the intersections of $\S^1\times\SU(2)$ with all
conjugates of $\S^1_{s,t}$. Since obviously both singular isotropy
groups are 1-dimensional they are $\S^1_{s,t}$ and $\S^1_{-s,t}$,
and thus both $B_{\pm}$ have fundamental group $\Z_{|s|}$.

We now combine the information $|p_i + p_j| \le 2$ and $|p_i+p_k| =
|p_k + p_j|$ gained so far, with the conditions for positive
curvature. One easily sees that this is only possible if either
$\bar{p} = (1,1,p)$ with $p>0$ or $\bar{p} = (-k,1,k-2)$ with $k$ odd
and $k\ge 5$. To exclude the latter case we use the fact that by
\cite{gwz} the only manifolds where $\S^1 \times \SU(2)\times \SU(2)$
acts by \co one are the Eschenburg spaces $E_p$ for some $p$. We will
now see that this contradicts topological invariants for these
spaces. Any Eschenburg space has finite, cyclic fourth cohomology
group, and for $E_{\bar{p}} \in \mathcal{E}_2$, this order is
$r=|p_1p_2+p_1p_3+p_2p_3|$. Kruggel also computed the first Pontrjagin
class of an Eschenburg space \cite{kruggel} and it follows in
particular that $p_1(E_{\bar{p}}) \equiv 2(p_1+p_2+p_3)^2 \mod r$. The
\com $E_p$ therefore has $r=2p+1$ and $p_1\equiv p+5 \mod r$ whereas
$E_{(-k,1,k-2)}$ satisfies $r= (k-1)^2 + 1$ and $p_1\equiv \frac{}{} 2
\mod r $. This yields the desired contradiction.
\end{proof}

\smallskip

In the remaining  cases where the $\T^2 \times \SU(2)$ action does
not extend to one of lower cohomogeneity, it easily follows that:

\begin{thm}[Cohomogeneity two Case]\label{isoE2}
The identity component of the isometry group of an $E_{\bar{p}} \in
\mathcal{E}_2 -(\mathcal{E}_1 \cup \mathcal{A})$ is given by
$$
   \Iso_0(E_{\bar{p}}) = \T^2 \times \S,
$$
where $\S = \SO(3)$ when all $p_i$ are odd, and $\S = \SU(2)$ otherwise.
\end{thm}

\smallskip

As in the previous cases, there exists another component of the
isometry group generated by complex conjugation in $\SU(3)$. We
suspect that this is then the full isometry group, but were not able
to prove it.

\bigskip

We end this section with a brief discussion of the isometry group of
Eschenburg's so-called \emph{twisted flag}. This is the $\T^2$
biquotient of $\SU(3)$ defined by
$$E^6 := \SU(3)/\!/\T^2, \text{where} \T^2 = \{ (\diag(z, w, zw),
\diag(1, 1, z^2 w^2)) \mid  z, w \in \U(1) \}.$$

\no The metric on $E^6$ induced from the $\Ad \U(2)$ invariant metric
on $\SU(3)$ used in the previous sections has positive curvature
(cf.\ \cite{esch2}). The right action of $\U(2) = \S(\U(2)\U(1))
\subset \SU(3)$ commutes with the $\T^2$ action and in particular
induces an isometric action on $E$. Moreover, as before, $E/\U(2) =
\T^2 \backslash \SU(3)/\U(2) = \T^2 \backslash \CP^2$ is a triangle
and thus $E$ has cohomogeneity two. There are in fact no more
isometries in the identity component of the isometry group.

\begin{prop}
The identity component $\Iso_{\subo}(E^6)$ of the twisted  flag $E^6$
is $\U(2)$.
\end{prop}

\begin{proof}
In \cite{shankar4} it was shown that $E$ is not homotopy equivalent to
any homogeneous space and in  \cite{searle}  that $E$ does not support
any positively curved cohomogeneity one metric. So any potential
extension of the $\U(2)$ action is again by cohomogeneity two.

From the rank rigidity theorem \eqref{rank} we know that $\Iso(E^6)$
has rank at most 2, and from the degree theorem \eqref{degree} that
$\dim \Iso(E^6) \le 6$. The only rank 2 group containing $\U(2)$ of
dimension at most 6 is $\SO(4)$. If the $\U(2)$  action were to extend
to $\SO(4)$, it would have to have the same orbits by the submetry
lemma. But the 4-dimensional principal $\U(2)$ orbits have infinite
fundamental group whereas a 4-dimensional quotient of $\SO(4)$ has
finite fundamental group.

Finally, since the commuting actions by $\U(2)$ and $\T^2$ on $\SU(3)$
only have $\id$ in common, the induced $\U(2)$ action on $E^6$ is
effective.
\end{proof}

\begin{rem*}
It is interesting to note that the normal homogeneous Aloff--Wallach
space, $E_{1,1,-2}$, admits actions of any possible cohomogeneity $k$,
with $0 \le k \le 7$. The ``generic" Eschenburg space
$E_{\bar{a},\bar{b}} \in \mathcal{E} -  \mathcal{E}_2$ has
cohomogeneity 4 with respect to the natural $\T^3$ action.
\end{rem*}

\section{Fundamental Groups}

It is a simple and well known fact that any finite group is the
fundamental group of a non-negatively curved manifold. A basic
question is if there are any obstructions in positive curvature other
than finiteness of the fundamental group.

In analogy with the situation in negative curvature, where a theorem
of Preismann asserts that any abelian subgroup must be cyclic, Chern
proposed the same obstruction for positive curvature in
\cite{chern}. However, in \cite{shankar}, it was shown that two well
known positively curved manifolds, the Aloff--Wallach space $A_{1,1}$
and the cohomogeneity one Eschenburg space $E_{1,1,2}$, both admit a
free, isometric $\SO(3)$-action. In particular, any finite group $\F
\subset \SO(3)$ containing $\Z_2 \oplus \Z_2 = \S(\O(1)\O(1)\O(1))$ is
the fundamental group of a positively curved manifold contradicting
Chern's conjecture.  Soon after infinitely many examples with
fundamental group $\Z_3\oplus \Z_3$ were found in \cite{grovshank}
(one of these was also found independently in \cite{bazaikin2}). These
as well as other groups were in fact the first non-\emph{space form
groups}, to be exhibited in positive curvature. We proceed to
exhibit an abundance of examples of positively curved manifolds with
non-space form groups as fundamental groups.

In all cases we have encountered, we only get something interesting
from subgroups of the identity component. Moreover, as observed in
\cite{grovshank}, if $\Gamma \subseteq \Iso_0$ acts freely, then
$\Gamma$ intersects any maximal torus in a cyclic group.  According to
Borel \cite{borel} a compact, connected group $\G$ has a non-toral
$\Z_p\oplus \Z_p$ if and only if $\pi_1(\G)$ has $p$-torsion. From our
description of isometry groups here, this already suggests not to
expect more interesting free, finite sub-actions from
$\U(2)\times\SO(3) = \Iso_0(E_p)$ than from $\SU(2)\times
\SO(3)$. Similarly, only when all $p_i$ are odd, is there a chance
that $\Iso_0(E_{\bar{p}}) = \T^2 \times \SO(3)$ has interesting free
finite sub-actions, and they all come already from $\S^1\times
\SO(3)$.

Since $\SO(3)$ has many interesting finite subgroups that are not
space form groups, it is worthwhile to determine for which
$E_{\bar{p}} \in \mathcal{E}_2$ an $\SO(3) \subset \Iso(E_{\bar{p}})$
acts freely. Moreover, $\SO(3)$ is a normal subgroup in
$\Iso(E_{\bar{p}})$, so the isotropy groups of the action are simply
its intersection with the regular and singular isotropy groups of the
\co two action. By considering the isotropy groups corresponding to
the vertex points determined in \eqref{lensspaces}, it follows that if
the action is free, then $\bar{p} = (1,1,1)$ or $\bar{p} =
(1,1,2)$. This recovers the main result in \cite{shankar} (see also
\cite{chaves} where it is shown that $A_{1,1}$ is a fat,
$\SO(3)$-principal bundle over $\CP^2$).

\begin{thm}
$\SO(3) \subseteq \Iso(E_{\bar{p}})$ acts freely if and only if
$\bar{p} = (1,1,1)$ or $\bar{p} = (1,1,2)$.
\end{thm}

\smallskip

Although $\SO(3)$ does not act freely on other Eschenburg spaces, we
will see that any finite subgroup $\F \subset \SO(3)$ in fact acts
freely on a large class of these spaces. As we will see, the same
assertion holds for the product of many cyclic groups with $\F \subset
\SO(3)$.

\bigskip

We begin with the cohomogeneity one spaces $E_p$, where
$\Iso(E_p)_{\subo} \cong \SO(3)\times \U(2)$. For any finite group
$\Gamma \subset \SU(2)\times \SU(2)$, we let $\Gamma(p) \subset
\Iso(E_p)$ be the image. Recall that $\H =
\{(\pm{\id})^{p+1},(\pm{\id})^{p}\}$ is the kernel of the action.

\begin{prop} For a finite group $\Gamma \subset \SU(2)\times \SU(2)$,
the group $\Gamma(p)$ acts freely on $E_p$ if  and only if for all
$(\gamma_1, \gamma_2) \in \Gamma - \H$ we have the conditions
\begin{itemize}
\item $|\gamma_1| \ne |\gamma_2|$, and $(\gamma_1, \gamma_2) \ne \pm
(-\id, \id)$, and

\item$|\gamma_1|$ does not divide $p|\gamma_2|$ or $|\gamma_2|$ does
not divide $(p+1)|\gamma_1|$,
\end{itemize}

\no where $|\gamma_i|$ denotes the order of $\gamma_i$.
\end{prop}

\begin{proof}
Any element of $\SU(2)$ has eigenvalues $\{\lambda, \bar{\lambda}\}$
for some $\lambda \in \U(1)$; let $\gamma_i$ have eigenvalues
$\{\lambda_i, \bar{\lambda}_i\}$. Clearly $\Gamma(p)$ fails to act
freely only if for some $(\gamma_1, \gamma_2) \in \Gamma - \H$ there
is a $(\diag(z,z,z^{p}), \diag(1,1,z^{p+2}) \in \S^1_p$ such that
$\diag(z,z,z^{p}) \gamma_1$ and $\diag(1,1,z^{p+2}) \gamma_2$ are
conjugate, i.e., have the same set of eigenvalues. So the action is
free unless the sets $\{z\lambda_1, z\bar{\lambda}_1, z^p\}$ and
$\{\lambda_2, \bar{\lambda}_2, z^{p+2}\}$ are the same. This happens
only if either: (a)   $ \lambda_1 = \pm {\lambda_2}$,  or  $\lambda_1 = \pm
\bar{\lambda}_2$  or  (b) $\lambda_1^p = \lambda_2^{p+1}$ or
$\lambda_1^p = \bar{\lambda}_2^{p+1}$.
Case (a) corresponds to   $|\gamma_1| = |\gamma_2|$ or $(\gamma_1,
\gamma_2) = \pm (-\id, \id)$ while in case (b) we have that
$|\gamma_2|$ divides $|\gamma_1|(p+1)$ and $|\gamma_1|$ divides
$|\gamma_2|p$.
\end{proof}

This has some interesting consequences.  For instance, if $\Gamma
\subset \SU(2) \times \{\id\}$, and $p$ is even, then $\Gamma(p)$ acts
freely if none of the orders of elements of $\Gamma\backslash
\{-\id\}$ divide $p$. Also, $\Gamma(p) \times \Z_q$ acts freely on
$E_p$, $p$ even if $\gcd(|\Gamma|,q)=\gcd(|\Gamma|,p)=1$ or
$\gcd(|\Gamma|,q)=\gcd(|\Gamma|,p+1)=1$ We deal with $\Gamma(p)$ and
$\Z_q \times \Gamma(p)$ for $\Gamma \subset \{\id\} \times \SU(2)$,
and $p$ is odd similarly. This proves the cohomogeneity one part of
Theorem B in the introduction. Note that the quaternion group $\Gamma
= \{\pm1,\pm i, \pm j, \pm k\} \subset \SU(2)$ corresponds to
$\Gamma(p) = \Z_2 \oplus \Z_2 \subset \SO(3)$, and that $\Z_q \times
(\Z_2 \oplus \Z_2) = \Z_{2q} \oplus \Z_2$.

\bigskip

For the cohomogeneity two Eschenburg spaces $E_{p_1, p_2, p_3}$ with
all $p_i$ odd, recall that the right action by $\SU(2)$ acts
effectively as $\SO(3)$. For these we have

\begin{prop}
Any finite group $\F \subset \SO(3)$ acts freely on an infinite family
of spaces $\E_{\bar{p}} \in \mathcal{E}_2$, in particular when all
$p_i$ are odd and distinct and $p_i \equiv 1 \mod |\F|$.
\end{prop}

\begin{proof}
As usual we let $\F^*$ be the inverse image of $\F$ by the map $\SU(2)
\to \SO(3) = \SU(2)/\{\pm \id\}$. In particular, the order of any
element of $\F^*$ divides the order of $\F$.

Let $p_1$, $p_2$, $p_3$ be three distinct primes, all congruent to 1
modulo $|\F|$ (infinitely many such primes exist by Dirichlet's
theorem on primes in arithmetic progression). As in the proof above we
see that the action of $\F$ on the space $E_{p_1, p_2, p_3}$ fails to
be free only if the sets $\{z^{p_1}, z^{p_2}, z^{p_3}\}$ and
$\{\gamma, \bar{\gamma}, z^{\sum p_i}\}$ are the same for some $z \in
\U(1)$.  Without loss of generality we may assume that $z^{p_1} =
\gamma$, $z^{p_2} = \bar{\gamma}$, $z^{p_1 + p_2} = 1$. Since in
particular $p_1$ and $p_2$ are relatively prime, we can find integers
$a$, $b$ such that $ap_1 + bp_2 = 1$. Then $z = z^{ap_1 + bp_2} =
\gamma^{a-b}$.  But since also $p_1$ and $p_2$ are both congruent to 1
modulo $|\F|$, we have: $\gamma = z^{p_1} = \gamma^{(a-b)p_1} =
\gamma^{a-b}$ and $\bar{\gamma} = z^{p_2} = \gamma^{(a-b)p_2} =
\gamma^{a-b}$. This implies $\gamma = \bar{\gamma} = \pm \id$ which is
simply the kernel of the action. Hence, $\F^*/ \{\pm \id\}= \F$ acts
freely and isometrically on each of these spaces.
\end{proof}

\providecommand{\bysame}{\leavevmode\hbox
to3em{\hrulefill}\thinspace}

\end{document}